\documentclass[12pt]{amsart}
\usepackage{geometry}                
\usepackage{graphicx}
\usepackage{amssymb}
\usepackage{epstopdf}
\usepackage{xcolor}
\title{On the Phragmen-Lindel\"of theorem in strips}
\author{Kevin Smith}
\newtheorem*{theorem*}{Theorem}
\newtheorem*{lemma*}{Lemma}
\newtheorem*{definition*}{Definition}
\newtheorem*{proposition*}{Proposition}

\newtheorem*{corollary*}{Corollary}
\newtheorem*{conjecture*}{Conjecture}

\newtheorem{theorem}{Theorem}

\numberwithin{equation}{section}

\begin{document}
\maketitle
\begin{abstract}
In this paper we study an $L^{p}$ analogue of Bohr's abscissae of summability for Dirichlet series. For polynomially bounded analytic functions in a strip with order function $\mu$, convexity of $1/\mu$ is equivalent to approximate concavity of the abscissae in $p$. If $\mu$ obeys a functional equation of the Selberg class type, this is equivalent to the Lindel\"of hypothesis if $\mu'(1/2)$ does not exist. Otherwise, $\mu$ is everywhere differentiable (therefore subconvex) with quadratic decay near one.
\end{abstract}

\section{Introduction}
The Phragmen-Lindel\"of theorem tells us that if $s=\sigma+it$ and $f(s)$ is analytic in a strip, say $a\leq  \sigma\leq 1$, $t\geq 1$, then the order function 
\begin{eqnarray}\label{mu}
\mu(\sigma)=\inf\{\mu\geq 0:f(\sigma+it)\ll t^{\mu}\}\nonumber
\end{eqnarray}
is convex and monotone in the strip if $f$ is polynomially bounded on the edges. So if $\mu(a)>0$ and $\mu(1)=0$ there is a $a<b\leq 1$ such that 

\begin{eqnarray}\label{sc}
\mu(\sigma)\leq \mu(a)\frac{b-\sigma}{b-a}\hspace{0.5cm}(a\leq \sigma\leq b)\hspace{0.5cm}\textrm{and}\hspace{0.5cm}\mu(\sigma)=0\hspace{0.5cm}(b\leq \sigma\leq 1).\nonumber
\end{eqnarray}
Many applications in analytic number theory depend on properties of $\mu$. For example, the functional equation of an $L$-function of degree $d\geq 1$ in the extended Selberg class $\mathcal{S}^{\sharp}$ (see Kaczorowski and Perelli \cite{KP}, also Conrey and Ghosh \cite{Conrey}, Dixit \cite{Dixit} and Murty \cite{Murty}) gives
\begin{eqnarray}\label{fe}
\mu(\sigma)=\mu(1-\sigma)+d\left(\frac{1}{2}-\sigma\right),\hspace{0.5cm}\hspace{0.5cm}\mu(\sigma)=0\hspace{0.5cm}(\sigma\geq 1)
\end{eqnarray}
and the Lindel\"of hypothesis is the conjecture that $b=1/2$. Yet the precise value of $\mu(\sigma)$  is not known for any $0<\sigma<1$ for any such function. \\

In this paper we consider the relationship between $\mu$ and the $L^{p}$ abscissae
\begin{eqnarray}\label{betadef}
\sigma_p(\alpha)=\inf\left\{\sigma\leq 1: \left(\int_{1}^{\infty}|f(\sigma+it)|^{p}\frac{dt}{t^{\alpha}}\right)^{1/p}<\infty  \right\}\hspace{0.5cm}(1<\alpha\leq \mu(a)p).\nonumber\\
\end{eqnarray}
Denoting by $\mu^{-1}$ the (convex and decreasing) inverse function
\begin{eqnarray}
\mu^{-1}:\left[0,\mu(a)\right]\rightarrow[a,b],
\nonumber
\end{eqnarray}
in general we have 
\begin{eqnarray}
\sigma_p(\alpha)\rightarrow\mu^{-1}(\alpha/p)\hspace{1cm}(p\rightarrow\infty)\nonumber
\nonumber
\end{eqnarray}
uniformly for $1<\alpha\leq \mu(a)p$ (quantitative estimates are given \S\ref{pth1}). If $f$ is analytic on the half plane $\sigma\geq a$ and represented by a convergent Dirichlet series for $\sigma>1$, then $\sigma_{2k}(\alpha)$ $(k\in\mathbb{N}, \alpha\geq 2)$ is an upper bound for Bohr's abscissa of summability 
\begin{eqnarray}\label{riesz}
&&\inf\left\{\sigma\leq 1 :f^k(s)=\lim_{x\rightarrow\infty} \sum_{n< x}\left(1-\frac{\log n}{\log x}\right)^{\alpha/2-1}\frac{c_{n,k}}{n^{s}}\right\}\nonumber\\
&=&\inf\left\{\sigma\leq 1: \int_{\sigma-i\infty}^{\sigma+i\infty}\frac{f^k(s)ds}{s^{\alpha/2}}\hspace{0.2cm}\textrm{exists}  \right\}
\end{eqnarray}
where  $c_{n,k}$ are the Dirichlet coefficients of $f^k$ (see also \cite{HardyRiesz}, \S 6.4). Bohr's studies of the abscissae (\ref{riesz}) in \cite{Bohr3,Bohr1,Bohr2} exhibit such $f$ for which the Lindel\"of hypothesis fails dramatically: for any $0<b\leq 1$ and any convex decreasing function $\mu$ with left derivative $\mu'_{-}(b)\leq -1$, there exists a Dirichlet series having this $\mu$ as its order function and abscissa of convergence $0$  (see \cite{Bohr2} and Kahane \cite{Kahane}). It is therefore of interest to know general conditions under which $\mu$ has properties that are consistent with the Lindel\"of hypothesis. One such condition is \emph{approximate concavity}, defined as follows.

\begin{definition*} The family of non-decreasing functions
\begin{eqnarray}\label{domaine}
\sigma_p(\alpha):\mathbb{R}_{\geq \alpha/\mu(a)}\rightarrow [a,b]\nonumber
\end{eqnarray}
is approximately concave if
\begin{eqnarray}\label{ac}
\phi \sigma_{q}(\alpha)+(1-\phi)\sigma_{r}(\alpha)\leq \sigma_p(\alpha)+o(1)\hspace{1cm}(p\rightarrow\infty)\nonumber
\end{eqnarray}
when $p=\phi q+(1-\phi)r$ and $0\leq \phi\leq 1$ uniformly for $1<\alpha\leq \mu(a)\min(q,r)$.
\end{definition*}

Our main result is the following.

\begin{theorem}\label{th1}
If $f$ is analytic in the strip $a\leq \sigma\leq 1$, $t\geq 1$,  $0<\mu(a)<\infty$ and $\mu(1)=0$, then $1/\mu$ is convex if and only if $\sigma_p(\alpha)$ is approximately concave.
\end{theorem}

Since $\mu$ is convex the right derivative $\mu'_{+}$ satisfies $\mu'_{+}\geq \mu'_{-}$ everywhere, and if $1/\mu$ is convex it satisfies $\mu'_{+}\leq \mu'_{-}$. Therefore $\mu$ is differentiable except possibly at $b$ if $\sigma_p(\alpha)$ is approximately concave. In view of Bohr's result recalled above, this is a strong condition. Thus if $f$ is analytic on the upper half plane we may let $a\rightarrow -\infty$ so $\mu$ is differentiable everywhere except possibly at $b$. If $\mu$ also satisfies (\ref{fe}), then $\mu$ is differentiable on $(-\infty,b)\cup(1-b,\infty)$, so if $b\neq 1/2$ then $\mu$ is everywhere differentiable. This gives an immediate corollary.

\begin{corollary*}\label{cor1}Let $f$ be analytic on the upper half plane. If $\mu$ satisfies the functional equation (\ref{fe}) and $\mu'(1/2)$ does not exist, then $f$ satisfies the Lindel\"of hypothesis if and only if $\sigma_p(\alpha)$ is approximately concave. 
\end{corollary*}

Otherwise, $\mu$ is everywhere differentiable and convexity of $1/\mu$ gives
\begin{eqnarray}
\mu(\sigma)\leq\frac{d}{2+4\sigma}=d\left(\frac{1}{2}-\sigma\right)+\frac{2d\sigma^2}{1+2\sigma}\hspace{1cm}(0\leq \sigma<b)\nonumber
\end{eqnarray}
so 
\begin{eqnarray}\label{fc}
\mu(\sigma)=\mu(1-\sigma)-d\left(\sigma-\frac{1}{2}\right)\leq 2d(1-\sigma)^2\hspace{1cm}(b< \sigma\leq 1).\nonumber
\end{eqnarray}
This quadratic decay near one may be compared with the exponent $3/2$ in Vinogradov-Korobov bounds, which leads to a conjecture.

\begin{conjecture*}Let $f\in \mathcal{S}^{\sharp}$ and suppose that $f(1+it)\ll (\log t)^A$. Then $\sigma_p(\alpha)$ is approximately concave and
\begin{eqnarray}\label{bs1}
f(s)\ll t^{2d(1-\sigma)^2}(\log t)^{A}\hspace{1cm}\left(1/2\leq \sigma\leq 1\right).
\end{eqnarray}
\end{conjecture*}

This conjecture implies that either $f$ satisfies the Lindel\"of hypothesis or $\mu$ is everywhere differentiable, so subconvexity $\mu(1/2)<d/4$ is implicit. Moreover, if $f$ has an Euler product and a classical zero-free region, then (\ref{bs1}) yields a zero-free region of the form
\begin{eqnarray}
\sigma>1-B(\log t)^{-1/2}(\log\log t)^{-A/2}\hspace{1cm}(t>C) \nonumber
\end{eqnarray}
in which the factor $(\log t)^{-1/2}$ offers an improvement over $(\log t)^{-2/3}$ in Vinogradov-Korobov zero-free regions (see Ford \cite{Ford} and Richert \cite{Richert}).   \\

Approximate concavity does not generally determine $b$ or $\sigma_p(\alpha)$, although the functional equation and subconvexity together imply the following.

\begin{theorem}\label{th2} Let $f\in \mathcal{S}^{\sharp}$, $k\in\mathbb{N}$ and suppose that $\mu(1/2)<d/4$. Then $f$ satisfies the Lindel\"of hypothesis if and only if
\begin{eqnarray}
\sigma_{2k}(\alpha)= b\left(1-\frac{\alpha-1}{d k}\right)\hspace{1cm}\left(k\geq \frac{\alpha-1}{d} \right).\nonumber
\end{eqnarray}
\end{theorem}

Taking $\alpha=2$ in Theorem \ref{th2}, it follows from a theorem of Helson \cite{Hel} that if an entire $f\in \mathcal{S}^{\sharp}$ satisfies the Lindel\"of hypothesis then the abscissa of convergence of the Dirichlet series of $f^k$ is $(1-1/dk)/2$, which is a conjecture of Rubinstein \cite{Rub}. \\

\newpage 
Lastly, we note that
\begin{eqnarray}
\sigma_p(\alpha)=\inf\left\{\sigma\leq 1: \log \|f_{\sigma}\|_{\alpha,p}<\infty  \right\}\nonumber
\end{eqnarray} 
in which
\begin{eqnarray}
  \|f_{\sigma}\|_{\alpha,p}^p= \int_{1}^{\infty}|f(\sigma+it)|^{p}\frac{dt}{t^{\alpha}}\hspace{1cm}(\alpha>1,p\geq \max(1,\alpha/\mu(a)))\nonumber
\end{eqnarray}
and by logarithmic convexity of $L^p$ norms that
\begin{eqnarray}\label{lc}
\log \|f_{\sigma}\|^p_{\alpha,p}\leq \phi \log \|f_{\sigma}\|^q_{\alpha,q}+(1-\phi)\log \|f_{\sigma}\|^r_{\alpha,r}
\end{eqnarray}
for $\sigma>\sigma_{\max(q,r)}(\alpha)$. In comparison, approximate concavity of the abscissae
\begin{eqnarray}\label{ac2}
\sigma_p(\alpha)\geq \phi \sigma_{q}(\alpha)+(1-\phi)\sigma_{r}(\alpha)+o(1)\nonumber
\end{eqnarray}
exhibits a rather special symmetry with (\ref{lc}).\\

Theorems \ref{th1} and \ref{th2} are proved in \S\ref{pth1} and \S\ref{pth2}, respectively.

\section{Proof of Theorem \ref{th1}}\label{pth1}
We begin by bounding the range of $\sigma_p(\alpha)$ in terms of $\mu^{-1}$. 

\begin{lemma*}\label{lem1}For $1<\alpha\leq \mu(a)p$ we have 
\begin{eqnarray}\label{lembound}
a\leq \mu^{-1}\left(\frac{\alpha}{p}\right)\leq \sigma_p(\alpha)\leq \mu^{-1}\left(\frac{\alpha-1}{p}\right)\leq b.
\end{eqnarray}
\end{lemma*}

\begin{proof}The upper bound is immediate because if $\alpha>1$ then for 
\begin{eqnarray}
\sigma =\mu^{-1}\left(\frac{\alpha-1}{p}\right)+\delta\hspace{1cm}0<\delta< c-\mu^{-1}\left(\frac{\alpha-1}{p}\right)\nonumber
\end{eqnarray}
 there is a $\gamma>0$ such that 
\begin{eqnarray}
|f(\sigma+it)|^{p}\leq C^{p}t^{p\mu\left(\mu^{-1}((\alpha-1)/p)+\delta\right)+p\epsilon}\leq C^{p}t^{\alpha-1-p(\gamma-\epsilon)}\hspace{1cm}(t\geq 1)\nonumber
\end{eqnarray}
for every $\epsilon>0$, in which case 
\begin{eqnarray}\label{integral}
\left(\int_{1}^{\infty}|f(\sigma+it)|^{p}\frac{dt}{t^{\alpha}}\right)^{1/p}
\end{eqnarray}
is bounded uniformly in $p$. \\

As for the lower bound, for every $\delta>0$ there is a $D>0$ and a sequence $t_n\rightarrow\infty$ such that 
\begin{eqnarray}\label{low}
|f(\sigma+it_n)|>  Dt_n^{\mu(\sigma)-\delta}.
\end{eqnarray}
Moreover, for  every $\gamma,\epsilon>0$ and $0\leq t-t_n\leq 1$, we have
\begin{eqnarray}\label{eb}
|f(\sigma+it_n)|-|f(\sigma+it)| &\leq &\left|\int_{t_n}^t f'(\sigma+iv)dv\right|\nonumber\\
&\leq&(t-t_n)\max_{t_n\leq v\leq t}|f'(\sigma+iv)|\nonumber\\
&\leq&\frac{t-t_n}{\epsilon}\max_{t_n\leq v\leq t}\max_{|s-\sigma-iv|= \epsilon}|f(s)|\nonumber\\
&\leq&\frac{C}{\epsilon}(t-t_n)t_n^{\mu(\sigma-\epsilon)+\gamma}
\end{eqnarray}
by Cauchy's inequality. Noting that (\ref{eb}) is $\leq |f(\sigma+it_n)|/2$ provided that
\begin{eqnarray}
t-t_n\leq \frac{\epsilon}{2C} t_n^{-\mu(\sigma-\epsilon)-\gamma}|f(\sigma+it_n)|, \nonumber
\end{eqnarray}
by (\ref{low}) we have $|f(\sigma+it)|\geq |f(\sigma+it_n)|/2$ when
\begin{eqnarray}
t-t_n\leq \frac{D\epsilon}{2C} t_n^{\mu(\sigma)-\mu(\sigma-\epsilon)-\gamma-\delta}. \nonumber
\end{eqnarray}
Thus the integral in (\ref{integral}) is
\begin{eqnarray}\label{lower}
&>& \sum_{n> 0}\int_{t_n}^{t_n+(D\epsilon/2C) t_n^{\mu(\sigma)-\mu(\sigma-\epsilon)-\gamma-\delta}} \frac{|f(\sigma+it)|^{p}dt}{|\sigma+it|^{\alpha}}\nonumber\\
&\gg_p&\sum_{n> 0}t_n^{\mu(\sigma)p-\alpha-p\delta}\int_{t_n}^{t_n+(D\epsilon/2C) t_n^{\mu(\sigma)-\mu(\sigma-\epsilon)-\gamma-\delta}} dt\nonumber\\
&\gg &\epsilon\sum_{n> 0}t_n^{\mu(\sigma)p+\mu(\sigma)-\mu(\sigma-\epsilon)-\alpha-\gamma-p\delta}\nonumber
\end{eqnarray}
in which the exponent of $t_n$ is positive for sufficiently small $\gamma,\delta,\epsilon>0$ when $\mu(\sigma)p>\alpha$. That is, when $\sigma<\mu^{-1}(\alpha/p)$. 
\end{proof}

Since $\mu(\sigma)$ $(a\leq \sigma\leq b)$ is continuous, decreasing and non-negative, $1/\mu$ is convex if and only if $\mu^{-1}(1/p)$ is concave on the domain $1/\mu(a)\leq p<\infty$. Let us first assume the latter, thus
\begin{eqnarray}\label{conc}
\phi \mu^{-1}\left(\frac{\delta}{q}\right)+(1-\phi)\mu^{-1}\left(\frac{\delta}{
r}\right)\leq \mu^{-1}\left(\frac{\delta}{p}\right)
\end{eqnarray}
for $p=\phi q+(1-\phi)r$, $0\leq \phi\leq 1$ and $1<\delta\leq \mu(a)\min(q,r)$. Set $q\leq r$. If $q$ is bounded then $\delta\leq \mu(a)q$ is bounded, in which case
\begin{eqnarray}
\sigma_{q}\left(\delta\right)\leq b \hspace{0.5cm}\textrm{and}\hspace{0.5cm}    \sigma_p(\delta)\sim    \sigma_{r}\left(\delta\right)\rightarrow \mu^{-1}(0)=b\hspace{0.5cm}(p\rightarrow\infty)\nonumber
\end{eqnarray}
therefore
\begin{eqnarray}\label{conv}
\phi \sigma_{q}\left(\delta\right)+(1-\phi)\sigma_{r}\left(\delta\right)\leq  \sigma_p(\delta)+o(1)\hspace{1cm}(p\rightarrow\infty)
\end{eqnarray}
holds trivially. If $q\rightarrow\infty$, using the upper bound in (\ref{lembound}) on the left hand side of (\ref{conc}) and the lower bound on the right gives 
\begin{eqnarray}\label{cono}
\phi \sigma_{q}\left(\delta+1\right)+(1-\phi)\sigma_{r}\left(\delta+1\right)
\leq \sigma_p(\delta).\nonumber
\end{eqnarray}
Since $ \sigma_{q}\left(\delta+1\right)\sim \sigma_{q}\left(\delta\right)$ as $q\rightarrow \infty$ by (\ref{lembound}), (\ref{conv}) again follows. So $\sigma_p(\alpha)$ is approximately concave. 
\\

Conversely, we assume (\ref{conv}) holds for $1<\delta\leq \mu(a)q$ so 
\begin{eqnarray}
\frac{\delta}{r} \leq \frac{\delta}{p} \leq \frac{\delta}{q}\leq \mu(a)\nonumber
\end{eqnarray}
and put $\delta=p\alpha$ in (\ref{conv}), in which case 
\begin{eqnarray}
\frac{p\alpha}{r} \leq \alpha \leq \frac{p\alpha}{q} \leq \mu(a).\nonumber
\end{eqnarray}
Making 
\begin{eqnarray}
\frac{p}{r}\rightarrow\frac{\beta}{\alpha}\hspace{0.5cm}\textrm{and}\hspace{0.5cm}\frac{p}{q}\rightarrow\frac{\gamma}{\alpha}\hspace{0.5cm}(p,q,r\rightarrow\infty),\nonumber
\end{eqnarray}
the inequalities (\ref{lembound}) and (\ref{conv}) give
\begin{eqnarray}
\phi \mu^{-1}\left(\frac{p\alpha}{q}\right)+(1-\phi)\mu^{-1}\left(\frac{p\alpha}{r}\right)
\leq \mu^{-1}\left(\alpha\right)+o(1)\hspace{0.5cm}(p,q,r\rightarrow\infty)\nonumber
\end{eqnarray}
in which 
\begin{eqnarray}
\phi=\frac{1-p/r}{1-q/r}\rightarrow\frac{1-\beta/\alpha}{1-\beta/\gamma}.\nonumber
\end{eqnarray}
Thus
\begin{eqnarray}\label{limc}
\frac{1-\beta/\alpha}{1-\beta/\gamma}\mu^{-1}\left(\gamma\right)+
\frac{1-\gamma/\alpha}{1-\gamma/\beta}\mu^{-1}\left(\beta\right)
\leq \mu^{-1}\left(\alpha\right)\nonumber
\end{eqnarray}
where
\begin{eqnarray}\label{limc}
1/\alpha=\frac{1-\beta/\alpha}{1-\beta/\gamma}1/\gamma+
\frac{1-\gamma/\alpha}{1-\gamma/\beta}1/\beta,
\nonumber
\end{eqnarray}
so $\mu^{-1}(1/p)$ is concave on the domain $1/\mu(a)\leq p<\infty$.

\section{Proof of Theorem \ref{th2}}\label{pth2}

If $f$ satisfies the Lindel\"of hypothesis then $b=1/2$. An argument identical to (\cite{Titchmarsh}, \S 13.3) and a mean value theorem for Dirichlet polynomials (\cite{Montgomery}, \S 6) give 
\begin{eqnarray}
\lim_{T\rightarrow\infty}\frac{1}{T}\int_{1}^{T}|f(\sigma+it)|^{2k}dt= \sum_{n=1}^{\infty}\frac{|c_{n,k}|^2}{n^{2\sigma}} \hspace{1cm}(\sigma> 1/2)\nonumber
\end{eqnarray}
so
\begin{eqnarray}\label{lowerbound}
\textrm{meas}\{t\in[1,T]:|f(\sigma+it)|>\delta\}\gg T\hspace{1cm}(\sigma> 1/2)
\end{eqnarray}
 for all sufficiently small $\delta>0$. Using the functional equation we have
\begin{eqnarray}\label{functionalequation}
f(\sigma+it)\asymp t^{d(1/2-\sigma)}|f(1-\sigma-it)|\hspace{1cm}(t\rightarrow\infty)\nonumber
\end{eqnarray}  
so the lower bound (\ref{lowerbound}) implies that 
\begin{eqnarray}\label{lowerbound2}
\textrm{meas}\{t\in[1,T]:|f(\sigma+it)|\gg t^{d(1/2-\sigma)}\}\gg T\hspace{1cm}(\sigma< 1/2).\nonumber
\end{eqnarray}
It follows that the upper bound in (\ref{lembound}) is sharp, so
\begin{eqnarray}
\sigma_{2k}(\alpha)=\mu^{-1}\left(\frac{\alpha-1}{2k}\right)=\frac{1}{2}\left(1-\frac{\alpha-1}{d k}\right)\hspace{1cm}\left(k\geq  \frac{\alpha-1}{d}\right).\nonumber
\end{eqnarray}\\

Conversely, if 
\begin{eqnarray}
\sigma_{2k}(\alpha)=b\left(1-\frac{\alpha-1}{d k}\right)\hspace{1cm}\left(k\geq \alpha/d\right)\nonumber
\end{eqnarray}
we set $\alpha=1+2k\mu(\sigma)$ so making $k\rightarrow\infty$ in (\ref{lembound}) gives $\mu(\sigma)=d(1-\sigma/b)/2$. So if $\mu(1/2)<d/4$ then $1/2\leq b<1$ and (\ref{fe}) gives
\begin{eqnarray}\label{hold}
\mu(\sigma)=d\left(\frac{1}{2}-\sigma\right)\hspace{1cm}(0\leq \sigma\leq 1-b). 
\end{eqnarray}
Since $\mu$ is linear for $0\leq \sigma\leq b$, (\ref{hold}) holds also for $1-b\leq \sigma\leq b$, so $b=1/2$.

\end{document}